\documentclass{amsart}

\title{The Tropical Totally Positive Grassmannian}

\author{David Speyer and Lauren Williams} 
\address{Department of Mathematics, University of California, Berkeley, and
Department of Mathematics, MIT, Cambridge}

\usepackage{pstricks, pst-node, amssymb}
\psset{unit=1pt, arrowsize=4pt, linewidth=.7pt}
\psset{linecolor=blue}
\newgray{grayish}{.90}
\newrgbcolor{embgreen}{0 .5 0}

\def\vblack(#1, #2)#3{\cnode*[linecolor=black](#1, #2){3}{#3}}
\def\vwhite(#1,#2)#3{\cnode[linecolor=black,fillcolor=white,fillstyle=solid](#1,#2){3}{#3}}
\countdef\x=23
\countdef\y=24
\countdef\z=25
\countdef\t=26

\def\tbox(#1,#2)#3{
\x=#1 \y=#2 
\multiply\x by 12 
\multiply\y by 12 
\z=\x \t=\y
\advance\z by 12 
\advance\t by 12 
\psline(\x,\y)(\x,\t)(\z,\t)(\z,\y)(\x,\y)
\advance\x by 6
\advance\y by 6 
\rput(\x,\y){{\bf #3}}}

\usepackage{amsmath, amsthm, amssymb, amsbsy}
\usepackage{amsfonts, latexsym, stmaryrd, amscd, xy}
\usepackage[mathscr]{eucal}
\usepackage{epsfig}
\newtheorem{theorem}{Theorem}[section]

\newtheorem{proposition}[theorem]{Proposition}
\newtheorem{lemma}[theorem]{Lemma}

\newtheorem{corollary}[theorem]{Corollary}
\newtheorem{definition}[theorem]{Definition}
\newtheorem{conjecture}[theorem]{Conjecture}

\usepackage{amsfonts}
\usepackage{xy}
\usepackage{amssymb}
\usepackage{amsmath}

\newcommand{\Q}{\mathbb Q}
\newcommand{\R}{\mathbb R}
\newcommand{\RR}{\mathcal R}
\newcommand{\CC}{\mathcal C}
\newcommand{\T}{\mathcal T}

\newcommand{\C}{\mathbb C}
\newcommand{\VV}{\mathcal V}
\newcommand{\PPP}{\mathcal P}
\newcommand{\PP}{\mathbb{P}}

\newcommand{\AAA}{\mathcal{A}}

\newcommand{\congruent}{\equiv}

\DeclareMathOperator{\Trop}{Trop}
\DeclareMathOperator{\Path}{Path}
\DeclareMathOperator{\Web}{Web}
\DeclareMathOperator{\Sum}{Sum}
\DeclareMathOperator{\Prod}{Prod}

\DeclareMathOperator{\val}{val}

\DeclareMathOperator{\init}{in}
\DeclareMathOperator{\Spec}{Spec}
\DeclareMathOperator{\id}{id}

\newcommand{\thmrefer}[1]{\renewcommand\thetheorem
  {\protect\ref{#1}}\addtocounter{theorem}{-1}}

\xyoption{all}
\CompileMatrices

\begin{document}

\begin{abstract}
Tropical algebraic geometry is the geometry of the tropical semiring
$(\R, \min, +)$.  The theory of total positivity
is a natural generalization of the
study of matrices with all minors positive. 
In this paper we introduce the totally positive part of the 
tropicalization of an arbitrary affine variety, 
an object which has the structure of a polyhedral fan.  
We then investigate the case of the Grassmannian, 
denoting the resulting fan $\Trop^+ Gr_{k,n}$.
We show that $\Trop^+ Gr_{2,n}$ is the Stanley-Pitman fan, which is
combinatorially the fan dual to 
the (type $A_{n-3}$) associahedron, and that $\Trop^+ Gr_{3,6}$ and 
$\Trop^+ Gr_{3,7}$ are closely related to the fans dual to the types 
$D_4$ and $E_6$ associahedra.  These results are strikingly reminiscent
of the results of Fomin and Zelevinsky, and Scott, who showed that the
Grassmannian has a natural cluster algebra structure which is 
of types $A_{n-3}$, $D_4$, and $E_6$ for
$Gr_{2,n}$, $Gr_{3,6}$, and $Gr_{3,7}$. We suggest a general conjecture about the positive part of the tropicalization of a cluster algebra.
\end{abstract}

\maketitle

\section{Introduction}

Tropical algebraic geometry is the geometry of the tropical semiring
$(\R, \min, +)$.  Its objects are polyhedral cell complexes which behave
like complex algebraic varieties.  Although this is a very new field in 
which many basic questions have not yet been addressed (see 
\cite{FirstSteps} for a nice introduction), tropical geometry has already
been shown to have remarkable applications to enumerative geometry (see \cite{Mikhalkin}), as well as connections to
representation theory (see \cite{Cluster1}, \cite{Cluster2}, \cite{Lusztig3}).

The classical theory of total positivity concerns matrices in which all minors
are positive.  However, in the past decade this theory has been 
extended by Lusztig (see \cite{Lusztig1} and \cite{Lusztig2}), who
introduced the totally positive variety $G_{>0}$ in an arbitrary reductive
group $G$ and the totally positive part $B_{>0}$ of a real flag variety
$B$.  In the process, Lusztig discovered surprising connections between his 
theory of canonical bases for quantum groups and the theory of total 
positivity.

In this paper we introduce the totally positive part (or positive part, 
for short) of the tropicalization of an arbitrary affine variety over
the ring of Puiseux series, and then investigate what we get in the case
of the Grassmannian $Gr_{k,n}$.  First we give a parameterization of the 
totally positive part of the Grassmannian, largely based on work of
Postnikov \cite{Postnikov}, and then we compute its tropicalization, 
which we denote by 
$\Trop^+ Gr_{k,n}$.  We identify 
$\Trop^+ Gr_{k,n}$  with a polyhedral subcomplex of the 
${n \choose k}$-dimensional Gr{\"o}bner fan of the ideal of Pl{\"u}cker relations,
and then show that this fan, modulo its $n$-dimensional lineality space, is combinatorially equivalent to an
$(n-k-1)(k-1)$-dimensional fan which we explicitly describe.  As a special
case, we show that $\Trop^+ Gr_{2,n}$ is  a fan which appeared in the 
work of Stanley and Pitman (see \cite{StanleyPitman}), which parameterizes
certain binary trees, and which is combinatorially equivalent to the 
(type $A_n$) associahedron.  We also show that $\Trop^+ Gr_{3,6}$ 
and $\Trop^+ Gr_{3,7}$ are  fans
which are closely related to the fans of the types $D_4$ and $E_6$ 
associahedra,
which were first introduced in \cite{Y-systems}.  These results are 
strikingly reminiscent of the results of Fomin and Zelevinsky \cite{Cluster2}, 
and Scott \cite{Scott}, who showed 
that the Grassmannian has a natural cluster algebra structure which is of 
type $A_n$ for $Gr_{2,n}$, type $D_4$ for $Gr_{3,6}$, and type $E_6$ for
$Gr_{3,7}$. 
(Fomin and Zelevinsky proved the $Gr_{2,n}$ case and stated the other
results; Scott worked out the cluster algebra structure of all Grassmannians
in detail.)  Finally, we suggest a general conjecture about the 
positive part of the tropicalization of a cluster algebra.

\section{Definitions}

In this section we will define the tropicalization and positive part of the 
tropicalization of an arbitrary affine variety over the ring of Puiseux series. 
We will then describe the tropical varieties that will be of interest to us.

Let $\CC=\bigcup_{n=1}^{\infty} \C((t^{1/n}))$ and $\RR=\bigcup_{n=1}^{\infty} \R((t^{1/n}))$ be the fields of Puiseux series over $\C$ and $\R$.  Every 
Puiseux series $x(t)$ has a unique lowest term $a t^u$ where 
$a \in \C^{*}$ and $u\in \Q$.  Setting 
$\val(f) = u$, this defines the  
{\it valuation map}
$\val : (\CC^*)^n \rightarrow \Q^n, (x_1, \dots , x_n) \mapsto 
(\val(x_1), \dots , \val(x_n))$.  
We define $\RR^+$ to be $\{ x(t) \in \CC | \text{ the 
coefficient of the lowest term of }x(t) \text{ is real and positive} \}$. We will discuss the wisdom of this definition later;
for practically all purposes, the reader may think of $\CC$ as if it were $\C$ and of $\RR^+$ as if it were $\R^+$.

Let $I \subset \CC[x_1, \ldots, x_n]$ be an ideal. We define the 
tropicalization of $V(I)$, denoted $\Trop V(I)$, to be the closure of the
image under $\val$ of 
$V(I) \cap (\CC^*)^n$, where $V(I)$ is the variety
 of $I$. Similarly, we define 
the positive part of $\Trop V(I)$, which we will denote as $\Trop^+ V(I)$, to be the closure of the
image under $\val$ of $V(I) 
\cap (\RR^+)^n$. Note that $\Trop V$ and $\Trop^+ V$ are slight abuses of notation; they depend on the affine space in which $V$ is embedded and not solely on the variety $V$.

If $f \in \CC[x_1, \ldots, x_n] \setminus \{ 0 \}$, let the {\it initial form}
$\init (f) \in \C[x_1, \ldots, x_n]$ be 
defined as follows: write $f=t^a g$ for $a \in \Q$ chosen as large as possible 
such that all powers of $t$ in $g$ are nonnegative. Then $\init (f)$ is the 
polynomial obtained from $g$ by plugging in $t=0$. If $f=0$, 
we set $\init (f)=0$.
If $w=(w_1, \ldots, w_n) \in \R^n$ then $\init_w (f)$ is defined 
to be $\init (f(x_i t^{w_i}))$. If $I \subset \CC[x_1, \ldots, x_n]$ 
then $\init_w (I)$ is the ideal generated by $\init_w (f)$ for all $f \in I$. 
It was shown in \cite{SpeyerSturmfels} that $\Trop V(I)$ consists of the 
collection of $w$ for which $\init_w (I)$ contains no monomials. The 
essence of this proof was the following:

\begin{proposition} \cite{SpeyerSturmfels} \label{lift} If $w \in \Q^n$ and $\init_w (I)$ contains no monomial 
then $V(\init_w (I)) \cap (\C^*)^n$ is nonempty and 
any point $(a_1, \ldots, a_n)$ of this variety can be lifted to a 
point $(\tilde{a}_1, \ldots, \tilde{a}_n)\in V(I)$ with the leading 
term of $\tilde{a}_i$ equal to $a_i t^{w_i}$. 
\end{proposition}

We now prove a similar criterion to characterize the points in $\Trop^+ V(I)$.

\begin{proposition} \label{poschar}
A point $w=(w_1, \ldots, w_n)$ lies in $\Trop^+ V(I)$ if and only if $\init_w (I)$ does not contain any nonzero polynomials 
in $\R^+ [x_1,\dots , x_n]$.
\end{proposition}

In order to prove this proposition, we will need the following result of 
\cite{PosPoint}, which relies heavily on a result of \cite{Hand}.

\begin{proposition}\cite{PosPoint}\label{Einsedler}
An ideal $I$ of $\R [x_1, \dots , x_n]$ contains a nonzero element
of $\R^+ [x_1, \dots , x_n]$ if and only if 
$(\R^+)^n \cap V(\init_\eta (I)) = \emptyset$ for all $\eta \in \R^n$.
\end{proposition}

We are now ready to prove Proposition \ref{poschar}.

\begin{proof}
Define $T \subset \Q^n$ to be the image of $V(I) \cap (\RR^+)^n$ under $\val$. 
Let $U$ denote the subset of $\R^n$ consisting of those $w$ 
for which $\init_w (I)$ contains no polynomials with all positive terms. By definition,  $\Trop^+ V(I)$ is the closure of 
$T$ in $\R^n$. We want to show that the closure of $T$ is $U$.  

It is obvious that $T$ lies in $U$.  $U$ is closed, as the property that $\init_w (f)$ has only positive terms is open as $w$ varies. Thus, the closure of $T$ lies in $U$.

Conversely, suppose that $w \in U$. Then, by Proposition \ref{Einsedler}, for some $\eta \in \R^n$, $(\R^+)^n \cap V(\init_{\eta} (\init_w (I))) \neq \emptyset$. For $\epsilon >0$ sufficently small, we have $\init_{\eta} (\init_w (I))=\init_{\epsilon \eta + w} (I)$. Therefore we can find a sequence $w_1$, $w_2$, \dots approaching $w$ with $(\R^+)^n \cap V(\init_{w_i} (I)) \neq \emptyset$.

As $w$ varies, $\init_w (I)$ 
 takes on 
only finitely many values, and the subsets of $\R^n$ on which $\init_w (I)$ 
takes a 
specific value form the relative interiors of the faces of a complete rational 
complex known as the Gr\"obner complex (see \cite{Sturmfels2}).  These complexes are actually fans when
$I$ is defined over $\R$ (\cite{Sturmfels}). Therefore, we may perturb each $w_i$, while preserving $\init_{w_i} (I)$, in order to assume that the $w_i \in \Q^n$ and we still have $w_i \to w$. Then, by Proposition \ref{lift}, each $w_i \in T$, so $w$ is in the closure of $T$ as desired.
\end{proof}

\begin{corollary}\label{fancorollary}
$\Trop V(I)$ and $\Trop^+ V(I)$ are closed subcomplexes of the 
Gr\"obner complex. In particular, they are polyhedral complexes. If $I$
is defined over $\R$, then  
$\Trop V(I)$ and $\Trop^+ V(I)$ are closed subfans of the 
Gr\"obner fan.
\end{corollary}

One might wonder whether it would be better to modify the 
definition of $\RR^+$ to require that our power series lie in $\RR$. This 
definition, for 
example, is more similar to the appearance of the ring of formal powers series 
in \cite{Lusztig3}. One can show that in the case of the Grassmanian,
this difference is unimportant. 
Moreover, the definition used here has the advantage that it makes it easy to 
prove that the positive part of the tropicalization is a fan. 

Suppose $V(I) \subset \CC^m$ and $V(J) \subset \CC^n$ are varieties and we 
have a rational map $f : \CC^m \to \CC^n$ taking $V(I) \to V(J)$. 
Unfortunately, 
knowing $\val (x_i)$ for $1 \leq i \leq m$ does not in general 
determine $\val (f(x_1, \ldots, x_m))$, so we don't get a nice
map $\Trop V(I) \to \Trop V(J)$. 
However, 
suppose that $f$ takes the 
positive points of $V(I)$ surjectively onto the 
positive points of $V(J)$ and suppose that $f=(f_1, \ldots, f_n)$ is 
{\it subtraction-free}, 
that is, the formulas for the $f_i$'s are rational functions in the $x_i$'s whose numerators and denominators have positive coefficients. 
Define $\Trop f: \R^m \to \R^n$ by replacing every $\times$ in $f$ 
with a $+$, every $/$ with a $-$, every $+$ with a $\min$ and every constant $a$ with $\val(a)$.

\begin{proposition} \label{blah}
Suppose $V(I)\subset \CC^m$ and $V(J)\subset \CC^n$ are varieties.  Let
$f:\CC^m \rightarrow \CC^n$ be a subtraction-free rational map
taking $V(I)$ to $V(J)$ such that 
$V(I)\cap (\RR^+)^m$ surjects onto $V(J) \cap (\RR^+)^n$.  Then
$\Trop f$ takes $\Trop^+ V(I)$ surjectively onto 
$\Trop^+ V(J)$.
\end{proposition}

\begin{proof}
This follows immediately from the formulas $\val(x+y)=\min(\val(x), \val(y))$ and $\val(xy)=\val(x)+\val(y)$ for $x$ and $y \in \RR^+$.
\end{proof}

We now define the objects that we will study in this paper.
Fix $k$ and $n$, and let $N = {n \choose k}$.  
Fix a polynomial ring $S$ in $N$ variables with coefficients in a commutative
ring.
The  {\it Pl{\"u}cker 
ideal}
$I_{k,n}$ is the homogeneous prime ideal in $S$ consisting
of the algebraic relations 
(called 
{\it Pl{\"u}cker relations})
among the $k \times k$ minors of any
$k \times n$-matrix with entries in a commutative ring.

Classically, the {\it Grassmannian} $Gr_{k,n}$ is the projective variety 
in $\PP_{\C}^{N-1}$ defined by the 
ideal $I_{k,n}$ of Pl{\"u}cker relations. We write $Gr_{kn}(\CC)$ for the variety in $\PP^{N-1}_{\CC}$ defined by the same equations. Similarly, we write $Gr_{k,n}(\R)$ for the real points of the Grassmannian, $Gr_{k,n}(\R^+)$ for the real positive points, $Gr_{k,n}(\RR^+)$ for those points of $Gr_{k,n}(\CC)$ all of whose coordinates lie in $\RR^+$ and so on. We write $Gr_{k,n}(\C)$ when we want to emphasize that we are using the field $\C$,
and use $Gr_{k,n}$ when discussing results that hold with no essential modification for any field. The \emph{totally positive Grassmannian} is the set $Gr_{k,n}(\R^+)$.
    
An element of $Gr_{k,n}$ can be represented by a full rank $k \times n$
matrix $A$.  If $K \in {[n] \choose k}$ we define the
{\it Pl{\"u}cker coordinate} $\Delta_K(A)$ to be
the minor of $A$ corresponding to the columns of $A$ indexed by $K$.
We identify the element of the Grassmannian with the matrix $A$ 
and with its set of Pl{\"u}cker coordinates (which satisfy the
Pl{\"u}cker relations).

Our primary object of study is the  {\it tropical positive Grassmannian} 
$\Trop^+ Gr_{k,n}$, which is a fan, by Corollary \ref{fancorollary}.
As in \cite{SpeyerSturmfels}, this fan has an $n$-dimensional lineality
space.  
Let $\phi$ denote the map from 
$(\C^*)^n$ into $(\C^*)^{n \choose k}$ which sends 
$(a_1, \dots , a_n)$ to the ${n \choose k}$-vector whose 
$(i_1, \dots, i_k)$-coordinate is $a_{i_1} a_{i_2} \cdots a_{i_k}$.
We abuse notation by also using $\phi$ for the same map $(\CC^*)^n \to (\CC^*)^{n \choose k}$.
Let $\Trop \phi$ denote the corresponding linear map 
which sends $(a_1, \dots , a_n)$ to the ${n \choose k}$-vector whose 
$(i_1, \dots, i_k)$-coordinate is $a_{i_1} + a_{i_2}+ \cdots + a_{i_k}$.
The map $\Trop \phi$ is injective, and its image is the common lineality 
space of all cones in $\Trop Gr_{k,n}$.

\section{Parameterizing the totally positive Grassmannian $Gr_{k,n}(\R^+)$} \label{Parameterize}

In this section we explain two equivalent ways to parameterize 
$Gr_{k,n}(\R^+)$, as well as a way to parameterize 
$Gr_{k,n}(\R^+) / \phi ( (\R^+)^n)$.
 The first
method, due to Postnikov \cite{Postnikov}, uses a certain directed graph 
$\Web_{k,n}$ with 
variables associated to each of its $2k(n-k)$  edges.
The second method is closely related to the first and uses the same
graph, but this time 
variables are associated to each of its 
$k(n-k)$ regions.  This has the advantage of giving a bijection
between $(\R^+)^{k(n-k)}$ and $Gr_{k,n}(\R^+)$.  Finally, we use $\Web_{k,n}$ 
with variables labelling each of its 
$(k-1)(n-k-1)$ inner regions in order to give a bijective parameterization
of $Gr_{k,n}(\R^+) /\phi((\R^+)^n)$.

\begin{figure}[h]
\centerline{\epsfig{figure=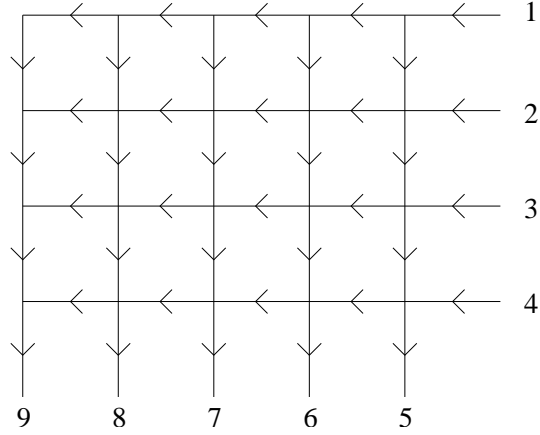}}
\caption{$\Web_{k,n}$ for $k=4$ and $n=9$}
\label{web49}
\end{figure}

Let $\Web_{k,n}$ be the directed graph which is obtained from a 
$k$ by $n-k$ grid, as shown in Figure \ref{web49}.  It has 
$k$ incoming edges on the right and $n-k$ outgoing edges on the bottom,
and the vertices attached to these edges are labelled clockwise
from $1$ to $n$.  We denote the set of $2k(n-k)$ edges by $E$.
Let us associate a formal variable
$x_e$ with each edge $e \in E$, and if there is no ambiguity,
we abbreviate the collection $\{x_e\}$ by $x$. If $p$ is a path
on $\Web_{k,n}$ (compatible with the directions of the edges),
then we let $\Prod_p(x)$ denote $\prod_{e \in p} x_e$.  And if $S$ is a set of 
paths on $\Web_{k,n}$, then we let $\Prod_S(x)$
denote $\prod_{p \in S} \Prod_p(x)$. 

As in \cite{Postnikov}, we define 
a $k \times n$ matrix $A_{k,n} (x)$, whose entries
$a_{ij}(x)$ are polynomials in the variables $x_e$, by the 
following equation:
\begin{equation*}
a_{ij}(x) = (-1)^{i+1} \sum_p \Prod_p(x), 
\end{equation*}
where the sum is over all directed
paths $p$ from vertex $i$ to vertex $j$.
Note that the $k \times k$ submatrix of $A_{k,n} (x)$ obtained by 
restricting to the first $k$ columns is the identity matrix.
In particular, $A_{k,n}(x)$ is a full rank matrix and hence
we can identify it with an element of $Gr_{k,n}$.
Also note that every element of the totally positive Grassmannian
$Gr_{k,n}(\R^+)$ has a unique matrix representative whose leftmost 
$k \times k$ submatrix is the identity.
We shall see that as the $\{x_e\}$ vary over
$(\R^+)^{2k(n-k)}$, the $A_{k,n}(x)$ range over all of $Gr_{k,n}(\R^+)$.

We now show it is possible to express the maximal minors (Pl{\"u}cker coordinates)
of $A_{k,n}(x)$ as
subtraction-free rational expressions in the $x_e$, as shown in 
\cite{Postnikov}.  If $K\in {[n] \choose k}$, then let 
$\Path(K)$ denote the set
\begin{equation*}
 \{S: S \text{ is a set of pairwise vertex-disjoint
paths  from }[k]\backslash (K\cap [k]) \text{ to }K\backslash ([k]\cap K) \}.
\end{equation*}
Note that for $K=[k]$, we consider the empty set to be a legitimate
set of pairwise vertex-disjoint paths.

Applied to $\Web_{k,n}$, 
Theorem 15.4 of \cite{Postnikov} implies 
the following.

\begin{proposition} \label{Plucker}
The Pl{\"u}cker coordinates of $A_{k,n}(x)$ are given by
\begin{equation*}
\Delta_K (A_{k,n}(x)) = \sum_{S \in \Path(K)} \Prod_S(x).
\end{equation*}
\end{proposition}

\begin{proof}
We give a brief proof of this result: the main idea
is to use the well-known Gessel-Viennot trick \cite{Gessel}.
First note that $a_{ij}(x)$ has a combinatorial interpretation: it is 
a generating function keeping track of paths from $i$ to $j$.
Thus, the determinant of a $k \times k$ submatrix of 
$A_{k,n}(x)$ corresponding to the column set $K$
also has a combinatorial interpretation: it is 
a generating function for all sets of paths from $[k]\backslash (K\cap [k])$ to 
$K\backslash ([k]\cap K)$,
with the sign of each term keeping track of the number of crossings
in the corresponding
path set.  What we need to show is that this is equal to the
sum of the contributions from path sets which are pairwise 
vertex-disjoint.
To see this, consider a path set which does have an intersection.
Look at its lexicographically last intersection,
and compare this path set to the one obtained from it by switching 
the two path tails starting at that point of intersection.  These two
path sets get different signs, but have equal weights, and hence
they cancel each other out.
\end{proof}

Proposition \ref{Plucker}
allows us to define a map $\Phi_0: (\R^+)^{2k(n-k)} \rightarrow Gr_{k,n}(\R^+)$ as follows.   
Let $K\in {[n] \choose k}$, and 
define $P_K: (\R^+)^{2k(n-k)}\rightarrow \R^+$ by
\begin{equation*}
P_K(x) := 
\sum_{S \in \Path(K)} \Prod_S(x).
\end{equation*}
Clearly if we substitute positive values for each $x_e$,
then 
$P_K (x)$ will be  positive.  We now define $\Phi_0$ by
\begin{equation*}
\Phi_0(x) = \{P_K(x)\}_{K \in {[n]\choose k}}. 
\end{equation*}
In other words, $\Phi_0$ is the map which 
sends a collection of positive real numbers
$\{x_e\}$ to the element of $Gr_{k,n}$ with Pl{\"u}cker coordinates
$P_K (x)$ (which is identified with the matrix $A_{k,n}(x)$).

By Theorem 19.1 of \cite{Postnikov}, the map $\Phi_0$ is 
actually surjective:  any point 
in $Gr_{k,n}(\R^+)$ can be represented as $A_{k,n}(x)$ for some positive choices
of $\{x_e\}$.  In summary, we have the following result
(which will also be a consequence of our Theorem \ref{Bij}).

\begin{proposition} \label{Parameter}
The map $\Phi_0: 
(\R^+)^{2k(n-k)} \rightarrow Gr_{k,n}(\R)$
is a surjection onto $Gr_{k,n}(\R^+)$. 
\end{proposition}

Unfortunately, the method we have just described
uses $2k(n-k)$ variables to parameterize
a space of dimension $k(n-k)$.  We will now explain how to do 
a substitution of variables which will reduce the 
number of variables to $k(n-k)$.

We define an {\it inner region} of $\Web_{k,n}$ to be a 
bounded component of the complement of 
$\Web_{k,n}$ (viewed as a subset of $\R^2$).  And we define an 
{\it outer region} of $\Web_{k,n}$ to be one of the extra inner
regions we would obtain if we were to connect vertices $i$ and $i+1$
by a straight line, for $i$ from $1$ to $n-1$.  A {\it region}
is an inner or outer region.  Note that there are
$k(n-k)$ regions, which we denote by $R$, and there are
$(k-1)(n-k-1)$ inner regions.

Let us label each  region $r\in R$  with a new variable $x_r$,
which we define to be the product of its counterclockwise edge variables 
divided by the product of its clockwise edge variables.

\begin{figure}[h]
\centerline{\epsfig{figure=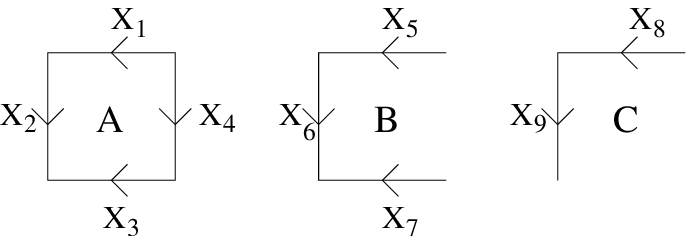}}
\caption{}
\label{Substitute}
\end{figure}

For example, the new variables $A, B, C$ shown in Figure \ref{Substitute}
would be defined by
\begin{equation*}
A:= \frac{x_1 x_2}{x_3 x_4}, \text{   }
B:= \frac{x_5 x_6}{x_7}, \text{   }
C:=  x_8 x_9.
\end{equation*}

It is easy to check that
for a path $p$ on $\Web_{k,n}$, 
$\Prod_p(x)$ is equal to the product of the variables attached to all
regions below $p$.  Since $\Prod_S(x)$ and $A_{k,n}(x)$ were defined  
in terms of the $\Prod_p(x)$'s, we can redefine these expressions in terms
of the $k(n-k)$ region variables.  
Proposition \ref{Parameter} still holds, but our map is now a map
$\Phi_1$ from $(\R^+)^{k(n-k)}$ onto $Gr_{k,n}(\R^+)$, taking the 
region variables $\{x_r\}$ to the element of $Gr_{k,n}(\R^+)$ represented by
$A_{k,n} (x)$.

Since we are now parameterizing a space of dimension $k(n-k)$
with $k(n-k)$ variables, we should have a bijection.  We shall prove 
that this is so by constructing the inverse map.

\begin{theorem} \label{Bij}
The map $\Phi_1: 
(\R^+)^{k(n-k)} \rightarrow Gr_{k,n}(\R^+)$, which maps
$\{x_r\}_{r\in R}$ to the Grassmannian element represented by $A_{k,n}(x)$,
is a bijection.
\end{theorem}

Before we prove this theorem, we need a lemma about 
matrices and their minors.  We use a very slight generalization
of a lemma which appeared in \cite{TP}.  For completeness,
we include the proof of this lemma. 
First we must define some
terminology.  Let $M$ be a $k \times n$ matrix.
Let $\Delta_{I,J}$ denote the minor of $M$ which uses row set $I$
and column set $J$.  We say that $\Delta_{I,J}$
is {\it solid} if $I$ and $J$ consist of several consecutive
indices; if furthermore $I\cup J$ contains $1$,
we say that $\Delta_{I,J}$ is {\it initial}.  Thus, an initial
minor is a solid minor which includes either the first column
or the first row. 

\begin{lemma} \cite{TP}\label{Minors}
A matrix $M$ is uniquely determined by its initial minors provided
that all these minors are nonzero.
\end{lemma}

\begin{proof}
Let us show that each matrix entry $x_{ij}$ of $M$ is uniquely 
determined by the initial minors.  If $i=1$ or $j=1$, there is 
nothing to prove, since $x_{ij}$ is an initial minor.  
Assume that $\min(i,j)>1$.  Let $\Delta$ be the initial
minor whose last row is $i$ and last column is $j$, and let
$\Delta^{\prime}$ be the initial minor obtained from $\Delta$
by deleting this row and column.  Then $\Delta = \Delta^{\prime} x_{ij}+P$,
where $P$ is a polynomial in the matrix entries $x_{i^{\prime} j^{\prime}}$
with $(i^{\prime}, j^{\prime}) \neq (i,j)$ and 
$i^{\prime} \leq i, j^{\prime} \leq j$.  Using induction on $i+j$,
we can assume that each $x_{i^{\prime} j^{\prime}}$ that occurs
in $P$ is uniquely determined by the initial minors, so the same is
true of $x_{ij} = (\Delta - P)/ \Delta^{\prime}$.   
\end{proof}

We now define a {\it reflected initial} minor to be a solid minor
$\Delta_{I,J}$ such that $I$ contains $k$ or $J$ contains
$1$.  Thus, a reflected initial minor is a solid minor which
includes either the first column or the last row.
A trivial corollary of Lemma \ref{Minors} is the following.

\begin{corollary}\label{silly}
A matrix $M$ is uniquely determined by its reflected
initial minors provided
that all these minors are nonzero.
\end{corollary}

Now we are ready to prove Theorem \ref{Bij}.

\begin{proof}
To prove the theorem, we will construct an explicit inverse map
$\Psi: Gr_{k,n}(\R^+) \rightarrow (\R^+)^{k(n-k)}$.  The first step
is to prove that $\Psi \Phi_1 = \id$.

Let us index the regions in $\Web_{k,n}$ by ordered pairs
$(i,j)$ as follows.
Given a region, we choose $i$ to be the label of the horizontal wire
which forms the upper boundary of the region, 
and choose $j$ to be the label of the
vertical wire which forms the left boundary of the region.
Now we define a map $K$ from the set of regions to ${[n] \choose k}$
by
\begin{equation*}
K(i,j):=\{1, 2, \dots , i-1\} \cup \{i+j-k, i+j-k+1, \dots , j-1, j \}.
\end{equation*}
If $(i,j)$ is not a region of $\Web_{k,n}$, then we
define $K(i,j):=\emptyset$.

Let $A$ be a $k \times n$ matrix whose initial $k\times k$ minor is the 
identity.  
We define $\Psi(A)$ by
\begin{equation}\label{preimage}
(\Psi (A))_{(i,j)}:=
\frac{ \Delta_{K(i,j)}(A) \Delta_{K(i+1,j-2)}(A) \Delta_{K(i+2,j-1)}(A) }
{ \Delta_{K(i,j-1)}(A)  \Delta_{K(i+1,j)}(A)  \Delta_{K(i+2,j-2)}(A) }.
\end{equation}
Note that by convention, we define $\Delta_{\emptyset}$ to be $1$.

See Figure \ref{inverse} for the definition of $\Psi$
in the case of 
$Gr_{3,6}(\R^+)$.  Note that for brevity, we have omitted the $A$'s from each term.

\begin{figure}[h]
\centerline{\epsfig{figure=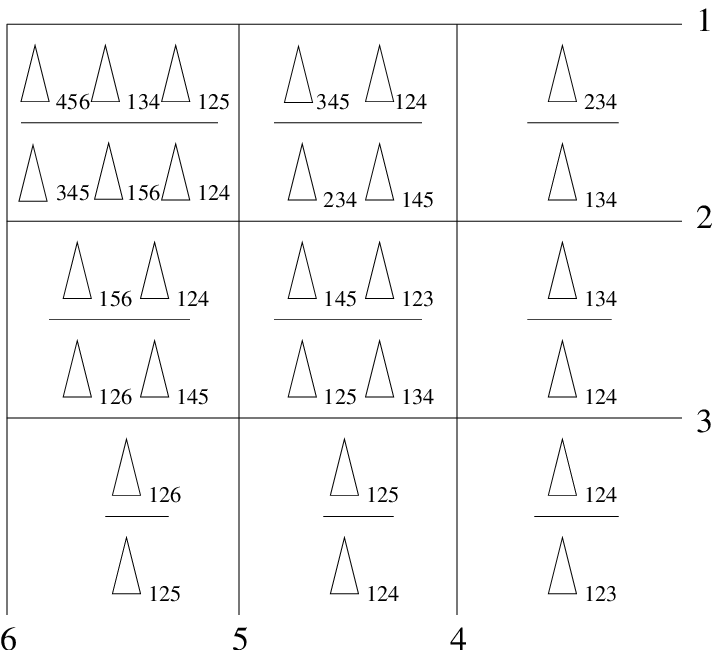}}
\caption{$\Web_{3,6}$}
\label{inverse}
\end{figure}

We claim that if  $A = \Phi_1 x$, then
$\Psi \Phi_1 x = x$.
To prove this, we note that the variable in region $(i,j)$ can
be expressed in terms of vertex-disjoint paths as follows.

First observe that if 
$K(i,j) \neq \emptyset$ then there is a unique set of pairwise vertex-disjoint
paths from $[k] \setminus ([k] \cap K(i,j))$ to 
$K(i,j) \setminus ([k] \cap K(i,j))$.  If one examines the terms in 
(\ref{preimage}) and draws in the six sets of 
pairwise vertex-disjoint paths on $\Web_{k,n}$ (say the three from the
numerator in red and the three from the denominator in blue) then it is
clear that every region in $\Web_{k,n}$ lies underneath an equal number
of red and blue paths -- except the region $(i,j)$, which lies underneath
only one red path.  Thus, by definition of the maps $P_K$,
it follows that
\begin{equation*}
(\Psi \Phi_1 (x))_{(i,j)} = \frac{ P_{K(i,j)}(x) P_{K(i+1,j-2)}(x) P_{K(i+2,j-1)}(x) }
{ P_{K(i,j-1)}(x)  P_{K(i+1,j)}(x)  P_{K(i+2,j-2)}(x) } = x_{(i,j)}.
\end{equation*}

%Note that since all of the values given by (\ref{preimage}) for 
%$A$ are positive numbers, this set of 
%region values is indeed the preimage of {\it some} point in $Gr_{k,n}(\R^+)$.

To complete the proof, it remains to show that $\Psi$ is injective.
This will complete the proof because we know that
$\Psi \Phi_1 \Psi = \Psi$, and $\Psi$ injective then implies that
$\Phi_1 \Psi = \id$.

Choose an element of $Gr_{k,n}(\R^+)$, which
we identify with its unique matrix representative $A$
whose 
leftmost $k \times k$ minor is the identity.
Let $\VV$ denote the set of rational expressions which
appear in the right-hand side of (\ref{preimage}) for all 
regions $(i,j)$ in $\Web_{k,n}$.  Let $\PPP$ denote the set of 
all individual Plucker coordinates which appear in $\VV$.
We prove that $\Psi$ is injective in two steps.  
First we show that the values of the expressions in $\VV$ 
uniquely determine the values of the
Pl\"ucker coordinates in $\PPP$.  Next we show that the values
of the Pl\"ucker coordinates in $\PPP$ uniquely determine 
the matrix $A$.

The first step is clear by inspection.  We illustrate the proof
in the case of $Gr_{3,6}(\R^+)$. 
By the choice of $A$,
$\Delta_{123}=1.$  Looking at the rational expressions in
Figure \ref{inverse}, we see that knowing the value 
$\frac{\Delta_{124}}{\Delta_{123}}$ 
determines $\Delta_{124}$; the value $\Delta_{124}$ together
with the value $\frac{\Delta_{134}}{\Delta_{124}}$ determines $\Delta_{134}$;
and similarly for $\Delta_{234}, \Delta_{125}, \Delta_{126}$.
Next, these values together with the value 
$\frac{\Delta_{145}\Delta_{123}}{\Delta_{125} \Delta_{134}}$
determines $\Delta_{145}$, and so on.  

For the second step of the proof, let 
$A^{\prime}$ denote the $k \times (n-k)$ matrix obtained from
$A$ by removing the leftmost $k \times k$ identity matrix.  
Note that the values of the Pl\"ucker coordinates
$\Delta_{K(i,j)}(A)$ (which are all elements of $\PPP$)
determine the values of all of the reflected initial minors of $A^{\prime}$.
(Each such Pl\"ucker coordinate is equal to one of the reflected initial minors,
up to sign.)  
Thus, by Corollary \ref{silly}, they uniquely determine the matrix
$A^{\prime}$ and hence $A$.
This completes the proof of Theorem \ref{Bij}.
\end{proof}

Now let us parameterize 
$Gr_{k,n}(\R^+) / \phi ( ({\R}^+)^n)$.
We shall show that we can
do this by using  variables corresponding to
only the $(k-1)(n-k-1)$  inner regions of 
$\Web_{k,n}$.  

First recall that the $n$-dimensional torus acts on $Gr_{k,n}(\R^+)$ by 
scaling columns of a matrix representative for $A\in Gr_{k,n}(\R^+)$.
(Although the torus has dimension $n$, this is actually just an
$(n-1)$-dimensional action as the scalars act trivially.)
Namely, 
\begin{equation*}
(\lambda_1, \dots , \lambda_n)    
\left(
\begin{matrix}
a_{11} & \dots & a_{1n} \\
\vdots && \vdots \\
a_{k1} & \dots & a_{kn}
\end{matrix} 
\right)
:=
\left(
\begin{matrix}
\lambda_1 a_{11} & \dots & \lambda_n a_{1n} \\
\vdots && \vdots\\
\lambda_1 a_{k1} & \dots & \lambda_n a_{kn}
\end{matrix} 
\right)
\end{equation*}
If $A\in Gr_{k,n}(\R^+)$ then we let $\overline{A}$ denote the torus orbit of $A$ under
this action.  
Note that if $K=\{i_1, \dots , i_k\}$, then 
$\Delta_K (\lambda A) = \lambda_{i_1} \lambda_{i_2} \dots \lambda_{i_k}
\Delta_K (A)$.

We will now determine the corresponding torus action on $(\R^+)^{k(n-k)}$ such that the above bijection commutes with the actions. If $r$ is an internal region then $x_r$ is a ratio of Pl\"ucker coordinates with the
same indices appearing on the top and bottom, so $x_r$ is not modified by the torus action. A simple 
computation shows that the torus acts transitively on the values of the outer region variables. Thus, taking the quotient by $\phi((\R^+)^n)$ on the right hand side of the equation corresponds to forgetting the outer variables on the left. 

Define a map $\Phi_2: (\R^+)^{(k-1)(n-k-1)} \rightarrow  Gr_{k,n}(\R^+) / \phi ( ({\R}^+)^n)$ by lifting a point $c \in (\R^+)^{(k-1)(n-k-1)}$ to any arbitrarily chosen point $\tilde{c} \in (\R^+)^{k(n-k)}$ and then mapping $c$ to $\overline{\Phi_1(\tilde{c})}$. We have just proven:

\begin{theorem} \label{SecondBijection}
The map 
$\Phi_2:(\R^+)^{(k-1)(n-k-1)}
\rightarrow  
Gr_{k,n}(\R^+) / \phi ( ({\R}^+)^n)$
is a bijection.
\end{theorem} 

\section{A fan associated to the tropical positive Grassmannian}\label{FAN}

In this section we will
construct a lower-dimensional
fan associated to the tropical positive  Grassmannian
$\Trop^+ Gr_{k,n}$. By methods precisely analogous to those above, we can prove an analogue of Theorem \ref{SecondBijection} for the field of Puiseux series.

\begin{theorem}
The map
$\Phi_2: (\RR^+)^{(k-1)(n-k-1)} \rightarrow
Gr_{k,n}(\RR^+) / \phi ( ({\RR}^+)^n)$
 is a bijection.
\end{theorem}

This theorem allows us to compute 
$\Trop^+ Gr_{k,n}/(\Trop \phi)(\R^n)$ by applying the valuation map to the image of $\Phi_2$. By Proposition \ref{blah}
we can tropicalize the map $\Phi_2$, obtaining the following surjective map.

\begin{equation*}
\Trop \Phi_2 : \R^{(k-1)(n-k-1)} \rightarrow \Trop^+ Gr_{k,n}/(\Trop \phi)(\R^n)
\end{equation*}

The map $\Trop \Phi_2$ is the map we get by replacing multiplication
with addition and addition with minimum in the definition
of $\Phi_2$.  Explicitly, it is defined as follows.
Let $K \in {[n] \choose k}$, and let inner region variables take on 
values $\{x_r\}$ in $\R$.  Outer region variables are  chosen arbitrarily.
If $p$ is a path on $\Web_{k,n}$
then
we let 
$\Sum_p(x)$ denote the sum of all variables which label
regions below $p$.  Similarly, if $S$ is a set of paths,
then let $\Sum_S(x)$ denote $\sum_{p\in S} \Sum_p(x)$.  
Now define $\Trop P_K (x) : \R^{(k-1)(n-k-1)} \rightarrow \R$ by
\begin{equation*}
\Trop P_K (x) := \min \{ \Sum_S(x): S\in \Path(K) \}.  
\end{equation*}
The map $\Trop \Phi_2$ is the map 
$$\Trop \Phi_2 :\R^{(k-1)(n-k-1)} \rightarrow \Trop^+ Gr_{k,n} / (\Trop \phi)(\R^n) \subset \R^N/(\Trop \phi)(\R^n)$$
given by
$$(\Trop \Phi_2 (x))_K=P_K(x).$$

\begin{definition}
The fan $F_{k,n}$ is the complete fan in $\R^{(k-1)(n-k-1)}$
whose maximal cones are the domains of linearity of the piecewise linear
map
$\Trop \Phi_2$.
\end{definition}

Because $\Trop \Phi_2$ surjects onto 
$\Trop^+ Gr_{k,n}/(\Trop \phi)(\R^n)$, the fan $F_{k,n}$ reflects the 
combinatorial structure of the fan $\Trop^+ Gr_{k,n}/(\Trop \phi)(\R^n)$,
which differs from $\Trop^+ Gr_{k,n}$ only through modding out by the linearity space.
However, $F_{k,n}$ is much easier to work with, as it lives
in $(k-1)(n-k-1)$-dimensional space as opposed to 
${n \choose k}$-dimensional space.

Since the maps $\Trop P_K$ are piecewise linear functions,
to each one we can associate a fan $F(P_K)$ whose 
maximal cones are the domains of linearity for $\Trop P_K$.
It is clear that the fan $F_{k,n}$ is the simultaneous refinement of 
all of the fans $F(P_K)$.

From now on, by abuse of notation, we will refer to 
$\Trop^+ Gr_{k,n}/(\Trop \phi)(\R^n)$ as $\Trop^+ Gr_{k,n}$.

\section{$\Trop^+ Gr_{2,n}$ and the associahedron}

In this section we will describe the fan $F_{2,n}$
associated to $\Trop^+ Gr_{2,n}$.
We show that this fan is exactly the Stanley-Pitman fan $F_{n-3}$, which
appeared in the
work of Stanley and Pitman in  
\cite{StanleyPitman}.
In particular, the face poset of $F_{2,n}$, with a top element $\hat{1}$
adjoined, is isomorphic to the face lattice of the normal fan of the
associahedron, a polytope whose vertices correspond
to triangulations of the convex $n$-gon. (In the language of \cite{CFZ}, this is the associahedron of type $A_{n-3}$.)

Let us first do the example of $\Trop^+ Gr_{2,5}$.
\begin{figure}[h]
\centerline{\epsfig{figure=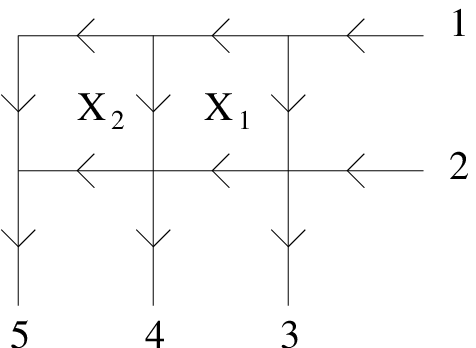}}
\caption{$\Web_{2,5}$}
\label{Web25}
\end{figure}

We use the web diagram $\Web_{2,5}$, as shown in Figure \ref{Web25}. 
The maps $\Trop P_K$ are given by:
\begin{align*}
\Trop P_{1j} &= 0 \text{ for all }j \\
\Trop P_{23} &= 0 \\
\Trop P_{24} &= \min(x_1,0) \\
\Trop P_{25} &= \min(x_1+x_2, x_1, 0) \\
\Trop P_{34} &= x_1 \\
\Trop P_{35} &= \min(x_1, x_1+x_2) \\
\Trop P_{45} &= x_1+x_2 \\
\end{align*}

\begin{figure}[h]
\centerline{\epsfig{figure=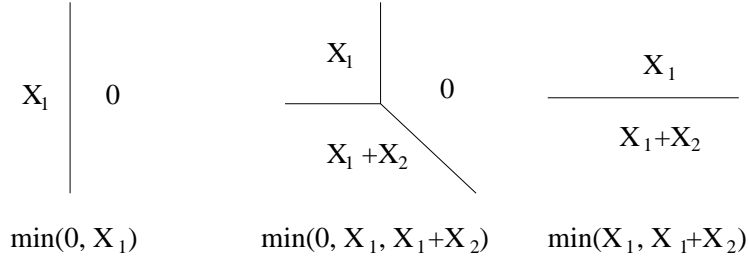}}
\caption{Fans for $\Trop P_J$}
\label{fan25}
\end{figure}

\begin{figure}[h]
\centerline{\epsfig{figure=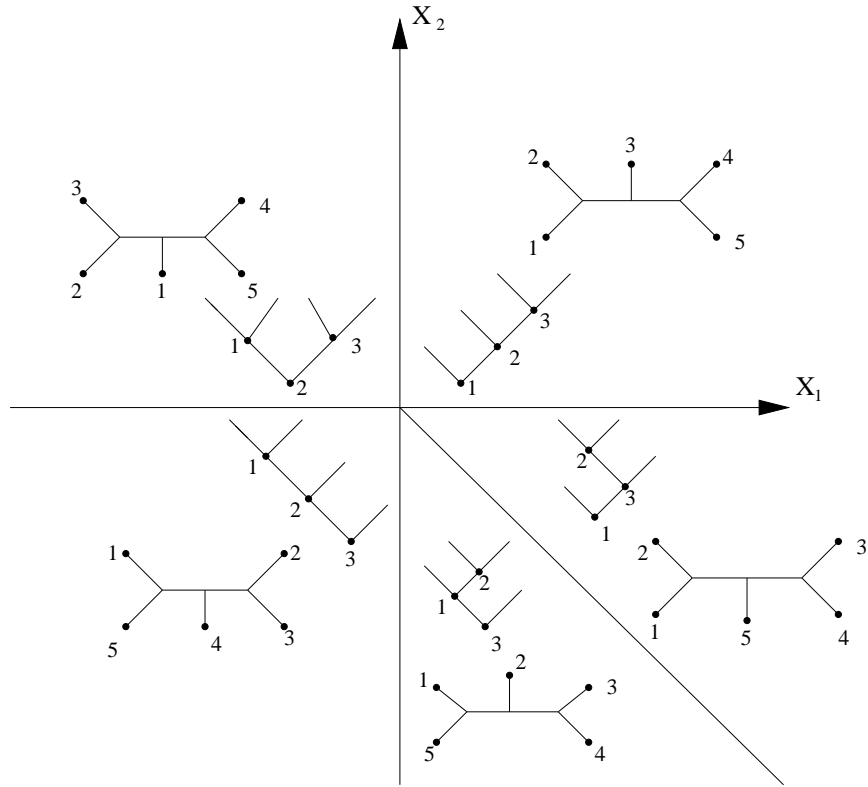}}
\caption{The fan of $\Trop^+ Gr_{2,5}$}
\label{BigFan25}
\end{figure}

Each map $\Trop P_K : \R^2 \rightarrow \R$ is piecewise linear
and so gives rise to the complete fan $F(P_K)$.  
For example, the map $\Trop P_{24}$ is linear on the 
region $\{ (x_1,x_2) : x_1 \geq 0 \}$,
where it is the function $(x_1,x_2) \mapsto 0$, and on the region
$\{(x_1,x_2): x_1\leq 0 \}$, where it is the function $(x_1,x_2) \mapsto x_1$.
Thus, $F(P_{24})$ is simply the subdivision of the real plane into
the regions $x_1 \geq 0$ and $x_1 \leq 0$.
The three nontrivial fans that we get from the maps $\Trop P_J$ are
shown in Figure \ref{fan25}.
In each picture,
the maximal cones of each fan are separated by solid lines.
$F_{2,5}$, which is the simultaneous refinement of the three nontrivial fans,
is shown in Figure \ref{BigFan25}.

In \cite{SpeyerSturmfels}, it was shown that maximal cones of the fan
$\Trop Gr_{2,n}$ correspond to trivalent trees on $n$ labelled leaves. 
It turns out that 
maximal cones of the fan $\Trop Gr_{2,n}^+$ correspond to 
trivalent planar trees 
on $n$ labelled leaves, as is illustrated in 
Figure \ref{BigFan25}.

We will now describe the fan  that appeared in \cite{StanleyPitman}, 
but first, we must review some notions about trees.
A {\it plane binary tree} is a rooted tree such that each vertex has either
two children designated as left and right, or none at all; and
an {\it internal vertex} of a binary tree
is a vertex which is not a leaf.  A {\it trivalent planar tree} is an 
(unrooted) tree
such that every vertex has degree three, and such that the leaves are labelled 
in a clockwise fashion.  It is known that both plane binary trees
with $n-1$ leaves, and trivalent planar trees with
$n$ labelled leaves, are counted by the Catalan number $c_{n-2}=\frac{1}{n-1} \binom{2(n-2)}{n-2}$.
 
There is a simple bijection between such  
trivalent planar trees and plane binary trees: 
if T a trivalent planar tree, then 
simply contract the edge 
whose leaf is labelled $1$, and make this the root.
This bijection is illustrated in Figure \ref{BigFan25}.

Let us now define the Stanley-Pitman fan $F_{n-3}$ in $\R^{n-3}$. 
(Note that we use different indices than are used in \cite{StanleyPitman}).
The maximal cones of $F_{n-3}$ are indexed by plane binary trees with 
$n-1$ leaves, in the following manner.
Let $T$ be a plane binary tree with $n-1$ leaves.
Label the internal vertices of $T$ with the numbers $1, 2,\dots n-2$
in the order of the first time we drop down to them from a child when doing a depth-first search from left to right starting
at the root. (See Figure \ref{BigFan25} for examples.) 
Let $x_1, \dots , x_{n-3}$ denote the coordinates in $\R^{n-3}$.
If the internal vertex $i$ of $T$ is the parent of vertex $j$, and
$i<j$, then associate with the pair $(i,j)$ the inequality
\begin{equation*}
x_{i} + \dots + x_{j-1} \geq 0,
\end{equation*}
while if $i>j$ then associate with $(i,j)$ the inequality
\begin{equation*}
x_{i} + \dots + x_{j-1} \leq 0.
\end{equation*}
These $n-3$ inequalities define a simplicial cone $C_T$ in 
$\R^{n-3}$.

The result proved in \cite{StanleyPitman} is the following.

\begin{theorem} [\cite{StanleyPitman}]
The $c_{n-2}$ cones $C_T$, as $T$ ranges over all plane binary trees
with $n-1$ leaves, form the chambers of a complete fan 
in $\R^{n-3}$.  Moreover, the face poset of $F_{n-3}$, with a top element 
$\hat{1}$ adjoined, is dual to the face lattice of the 
associahedron which parameterizes triangulations of the convex $n$-gon.
\end{theorem} 

The key step in proving that our fan $F_{2,n}$ is equal to the
Stanley-Pitman fan $F_{n-3}$ is the 
following lemma, also proved in \cite{StanleyPitman}.

\begin{lemma}[\cite{StanleyPitman}]
Let $D_i = \{(x_1, \dots , x_{n-3}) \in \R^{n-3} : x_1 + \dots + x_{i-1} =
  \min (0, x_1, x_1+x_2, \dots , x_1 + \dots + x_{n-3} ) \}$.
Let $\T_i$  consist of all plane binary trees with $n-1$ leaves and root $i$.
Then $D_i = \cup_{T \in \T_i} C_{T}$.
\end{lemma}

\begin{proposition}
The fan $F_{2,n}$  is equal to the fan $F_{n-3}$.
\end{proposition}

\begin{proof}
First let us describe the fan $F_{2,n}$ as explicitly as possible.
Note that if we label the regions of 
$\Web_{2,n}$ with the variables $x_1, \dots , x_{n-3}$ from right to 
left, then all of the maps
$\Trop P_K$ are of the form 
\begin{equation*}
\min (x_1+x_2+\dots + x_i, x_1+x_2+\dots + x_{i+1}, \dots , 
x_1+x_2+ \dots + x_j),
\end{equation*}
where $0 \leq i \leq j \leq n$.  Since this map 
has the same domains of linearity as
the map
\begin{equation*}
\theta_{ij}:=\min (x_i, x_i + x_{i+1}, \dots , 
x_i+ \dots + x_j),
\end{equation*}
 we can work with the maps $\theta_{ij}$ instead.
Let $F(i,j)$ be the fan whose cones are the domains of linearity of 
$\theta_{ij}$.  Then $F_{2,n}$ is the simultaneous refinement of 
all fans $F(i,j)$ where $0 \leq i \leq j \leq n$.

Now note that the previous lemma actually gives us an algorithm for
determining which cone $C_T$ a generic point $(x_1, \dots , x_{n-3}) \in 
\R^{n-3}$
lies in.
Namely, if we are given such a point, compute  the partial sums
of $x_1 + \dots + x_{i-1}$, for $1 \leq i \leq n-2$.  Choose $i$ such that 
$x_1 + \dots + x_{i-1}$ is the minimum of these sums.
(If $i=1$, the sum is $0$.)  Then the root of the tree
$T$ is $i$.  The left subtree of $T$ consists of vertices 
$\{1, \dots , i-1\}$, and the right subtree of $T$ consists
of vertices $\{i+1, \dots , n-2 \}$.  We now compute 
$\min \{0, x_1, x_1 + x_2, \dots , x_1 + \cdots + x_{i-2} \}$
and  
$\min \{x_{i}, x_{i}+x_{i+1}, \dots , x_{i} + \cdots + x_{n-3} \}$
in order to compute the roots of these two subtrees and so forth.

Now take a point $(x_1, \dots , x_{n-3})$ in a cone $C$ of $F_{2,n}$.
This means that the point is in a domain of linearity for all of the 
piecewise linear functions 
$\theta_{ij}=\min \{x_{i}, x_{i}+x_{i+1}, \dots , x_{i} + \cdots + x_j \}$
where $0 \leq i \leq j \leq n$, and we take $x_0$ to be $0$.
In other words, for each $i$ and $j$, there is a unique $k$ such that 
$x_i + \dots + x_k = 
\min \{x_{i}, x_{i}+x_{i+1}, \dots , x_{i} + \dots +x_j \}$.
In particular, we can reconstruct the tree $T$ such that 
$(x_1, \dots , x_{n-3}) \in C_T$, and {\it every} point 
$x \in C$ belongs to this same cone $C_T$.

Finally, we can show by induction that 
$C_T \subset C$.  (We need to show that all of the functions
$\theta_{ij}$ are actually linear on $C_T$.)
This shows that each cone $C$ in $F_{2,n}$ is actually equal to 
a cone $C_T$ in $F_{n-3}$, and conversely.

\end{proof}

\section{$\Trop^+ Gr_{3,6}$ and the type $D_4$ associahedron}

In connection with their work on cluster algebras, 
Fomin and Zelevinsky \cite{Y-systems} 
recently introduced certain polytopes
called {\it generalized associahedra} corresponding to
 each Dynkin type, of which the usual associahedron is the type
$A$ example.  When we computed $F_{3,6}$, the 
fan associated with
$\Trop^+ Gr_{3,6}$, we found that it was closely 
related to the normal fan of the type $D_4$ associahedron, in 
a way which we will now make precise.  (We defer the explanation 
of our computations to the end of this section.)

\begin{proposition}
The $f$-vector of $F_{3,6}$ is $(16, 66, 98,48)$.  The rays of 
$F_{3,6}$ are listed in Table \ref{16}, along with the 
inequalities defining the polytope that $F_{3,6}$ is normal to.
\end{proposition}

\begin{table}[h]
\begin{tabular}{|p{2.5cm}|p{3cm}|}
\hline
 $ e_1$& $x_1 \leq 5$\\
 $ e_2$& $x_2 \leq 7$ \\
 $ e_3$& $x_3 \leq 7$\\
 $ e_4$& $x_4 \leq 10$ \\
 $ -e_1$& $-x_1 \leq 0$\\
 $ -e_2$& $-x_2 \leq -2$ \\
 $ -e_3$& $-x_3 \leq -2$ \\
 $ -e_4$ & $-x_4 \leq -5$ \\
 $e_1 - e_2$& $x_1 - x_2 \leq 0$\\
 $e_1 - e_3$& $x_1 - x_3 \leq 0$ \\
 $e_1 - e_4$& $x_1 - x_4 \leq -1$ \\
 $-e_1 + e_4$& $-x_1 + x_4 \leq 9$ \\
 $e_2 - e_4$& $x_2 - x_4 \leq 0$ \\
 $e_3 - e_4$& $x_3 - x_4 \leq 0$\\
 $e_1 - e_2 - e_3$&$x_1 -x_2 -x_3 \leq -3$ \\
 $e_2+e_3-e_4$& $x_2+x_3-x_4 \leq 6$ \\
\hline
\end{tabular}
\bigskip
\caption{Rays and Inequalities for $F_{3,6}$}
\label{16}
\end{table}

Using the formulas of \cite{Y-systems}, we calculated the 
$f$-vector of the normal fan to the type $D_4$ associahedron: it
is $(16, 66, 100, 50)$.
More specifically, our fan has two cones which are of the form of a 
cone over a bipyramid. (Type FFFGG in the language of \cite{SpeyerSturmfels}.) 
If we subdivide these two bipyramids into two tetrahedra each, 
then we get precisely the $D_4$ associahedron. 

In Section \ref{cluster}, we will give some background
on cluster algebras and formulate a conjecture which explains the 
relation of $F_{3,6}$ to the normal fan to the type $D_4$ associahedron.

\begin{figure}[h]
\centerline{\epsfig{figure=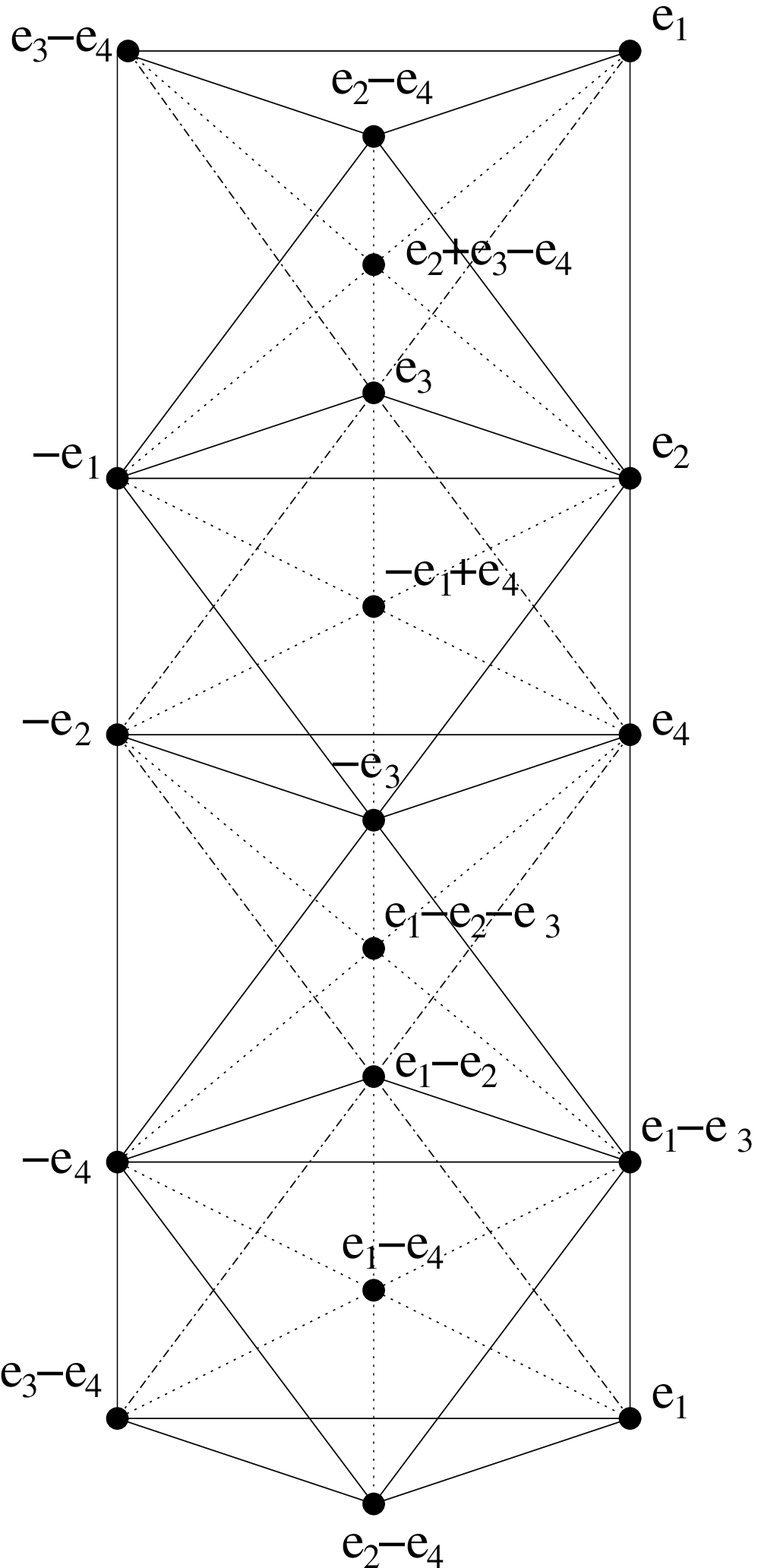}}
\caption{Glue bottom and top together to form a torus}
\label{D4_1}
\end{figure}

\begin{figure}[h] 
\centerline{\scalebox{0.7}{\epsfig{figure=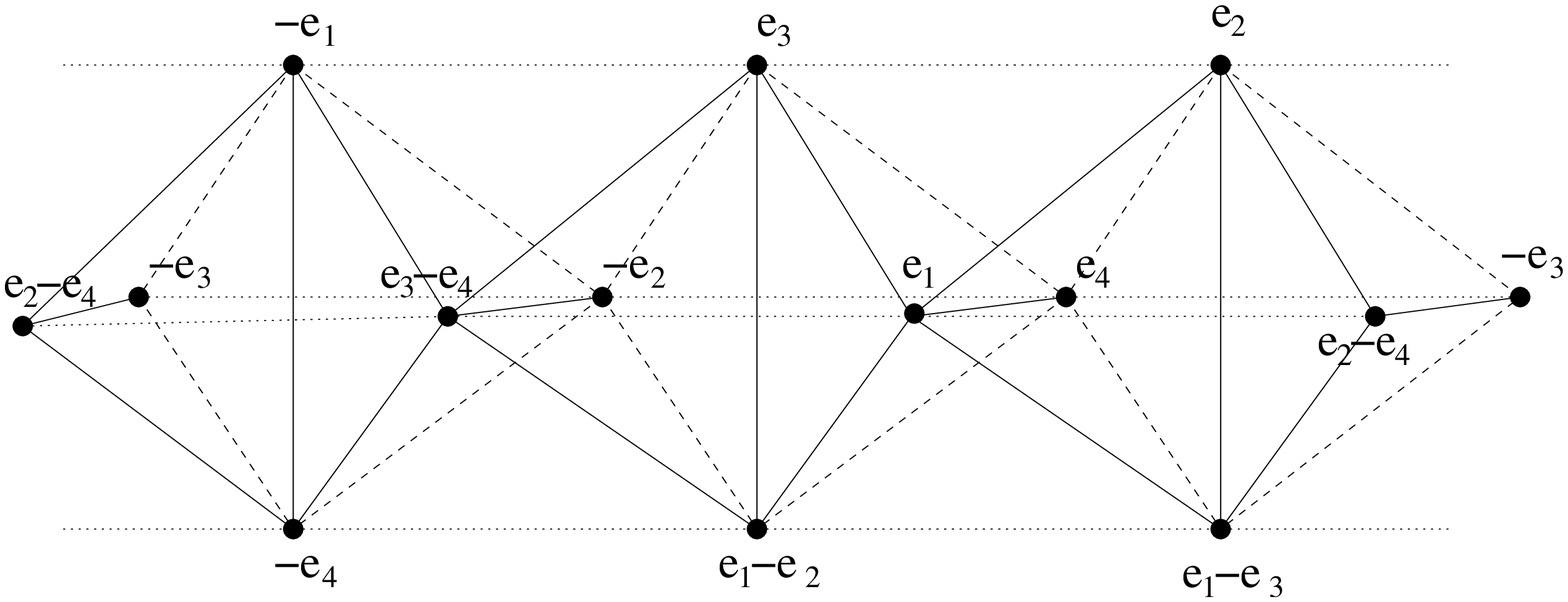}}}
\caption{Glue left and right ends together to form a torus}
\label{D4_2}
\end{figure}

We depict the intersection of $F_{3,6}$ with a 
sphere in Figures \ref{D4_1} and \ref{D4_2}. Each of the figures 
is homeomorphic to a solid torus, and the two figures glue together
to form the sphere $S^3$.  
The bipyramids in question have vertices $\{e_2+e_3-e_4, -e_1, e_2, e_3, -e_1+e_4\}$ and $\{e_1-e_2-e_3,-e_4,e_1-e_2,e_1-e_3, e_1-e_4\}$.

Now we will explain how we computed $F_{3,6}$.  We used two methods:
the first method was to use 
computer software (we used both cdd+ and Polymake) to compute the fan
which we described in Section \ref{FAN}.  The second method was to 
figure out which subfan of 
$\Trop Gr_{3,6}$ (which was explicitly described in \cite{SpeyerSturmfels})
was positive.  

To implement our first method, we used the 
well-known result that if $F_1$ and $F_2$ are fans which are normal
to polytopes $Q_1$ and $Q_2$, then the fan which is the 
refinement of $F_1$ and $F_2$ is normal to the Minkowski sum of 
$Q_1$ and $Q_2$.  Since the fan 
$F_{k,n}$
is the simultaneous refinement of all the fans $F(P_K)$, we 
found explicit coordinates for
 polytopes $Q(P_K)$ whose normal fans were the fans $F(P_K)$, and 
had the programs cdd+ and Polymake compute the Minkowski sum
$Q_{k,n}$ of all of
these polytopes. We then got explicit coordinates for the fan which
was normal to the resulting polytope.  

For the second method, we used the results in \cite{SpeyerSturmfels}:
we checked which of the  rays of $\Trop Gr_{3,6}$ did {\it not}
lie $\Trop^+ Gr_{3,6}$, and checked which facets of $\Trop Gr_{3,6}$
{\it did} lie in $\Trop^+ Gr_{3,6}$. As $\Trop^+ Gr_{3,6}$ is a closed subfan of $\Trop Gr_{3,6}$, this implied that every face of $\Trop Gr_{3,6}$ which lay in a totally positive facet was in $\Trop^+ Gr_{3,6}$ and every face of $\Trop Gr_{3,6}$ which contained a non-totally positive ray was not in $\Trop^+ Gr_{3,6}$; for every face of $\Trop Gr_{3,6}$, this proved sufficient to determine whether it was in $\Trop^+ Gr_{3,6}$ or not. 

%Although we will not write down all of the details here, 
%in Table \ref{Rays} we list certain rays of $\Trop Gr_{3,6}$ which are
%not totally positive, along with the Pl{\"u}cker relations which 
%demonstrate this.

%\begin{table}[h]
%\begin{tabular}{|p{3cm}||p{5cm}|}
%\hline
%$f_{1235}$ & $p_{134}p_{156}+p_{145}p_{136}=p_{135}p_{146}$\\
%$f_{1245}$ & $p_{123}p_{146}+p_{126}p_{134}=p_{124}p_{136}$\\
%$e_{124}$ & $p_{123}p_{245}+p_{125}p_{234}=p_{124}p_{235}$\\
%$e_{135}$ & $=p_{135}p_{245}$\\
%$g_{142536}$ & $p_{123}p_{246}+p_{126}p_{234}=p_{124}p_{236}$\\
%$g_{142635}$ & $Use 1435 for 4 of 5 terms$\\
%$g_{142356}$ & $p_{126}p_{456}+p_{156}p_{246}=p_{146}p_{256}$\\
%$g_{123546}$ & $Use 3456 and extra 1$\\
%\hline
%\end{tabular}
%\caption{Rays for $F_{3,6}$}
%\label{Rays}
%\end{table}

\section{$\Trop^+ Gr_{3,7}$ and the type $E_6$ associahedron}

As in the case of $F_{3,6}$, we used computer software to compute $F_{3,7}$,
the fan associated to 
$\Trop^+ Gr_{3,7}$. 

\begin{proposition}
The $f$-vector of $F_{3,7}$ is $(42, 392, 1463, 2583, 2163, 693)$.  Its rays  
are listed in Table \ref{42}, along with the inequalities defining the
polytope that $F_{3,7}$ is normal to.
Of the facets of this fan, $595$ are simplicial, $63$ have 7 vertices, $28$ have 8 vertices and $7$ have 9 vertices. All faces not of maximal dimension are simplicial.
\end{proposition}

Using the formulas of \cite{Y-systems}, we calculated 
the $f$-vector of the fan normal to the type $E_6$ associahedron: it is
$(42, 399, 1547, 2856, 2499, 833)$. 

In Section \ref{cluster}, we will explain why
$F_{3,7}$ differs from the $E_6$ fan, and how
one can refine  $F_{3,7}$ to get  a fan combinatorially equivalent
to the fan dual to the type $E_6$ associahedron. In this refinement, the simplicial facets remain facets. The $7$, $8$ and $9$ vertex facets split into $2$, $3$ and $4$ simplices respectively. The following table shows how the vertices of
the $7$, $8$, and $9$ vertex facets are grouped into simplices.

\begin{eqnarray*}
< ABCDEFG > & \Longrightarrow & < ABCDEF > \cup < ABCDEG > \\
< ABCDEFGH > & \Longrightarrow & < ABCDEF > \cup < ABCDFG > \\
 & & \phantom{< ABCDEF >} \cup <ABCDGH> \\
<ABCDEFGHI> & \Longrightarrow & <ABCDEF> \cup <ABCEFG> \cup \\ 
 & & <ABCFGH> \cup <ABCGHI> \phantom{\cup} 
\end{eqnarray*}

\begin{table}[h] 
\begin{tabular}{|p{3.5cm}|p{4cm}|}
\hline
 $ e_1 $ & $x_1 \leq 10$\\
 $ e_2 $ & $x_2 \leq 16$\\
 $ e_3 $ & $x_3 \leq 19$\\
 $ e_4 $ & $x_4 \leq 14$\\
 $ e_5 $ & $x_5 \leq 26$\\
 $ e_6 $ & $x_6 \leq 35$\\
 $ -e_1 $ & $-x_1 \leq -1$\\
 $ -e_2 $ & $-x_2 \leq -4$\\
 $ -e_3 $ & $-x_3 \leq -10$\\
 $ -e_4 $ & $-x_4 \leq -5$\\
 $ -e_5 $ & $-x_5 \leq -17$\\
 $ -e_6 $ & $-x_6 \leq -26$\\
 $e_1 - e_2$ & $x_1 - x_2 \leq -1$\\
 $e_1 - e_3$ & $x_1 - x_3 \leq -4$\\
 $e_1 - e_4$ & $x_1 - x_4 \leq -1$\\
 $e_1 - e_5$ & $x_1 - x_5 \leq -7$\\
 $e_1 - e_6$ & $x_1 - x_6 \leq -17$\\
 $e_2 - e_3$ & $x_2 - x_3 \leq -1$\\
 $e_2 - e_5$ & $x_2 - x_5 \leq -4$\\
 $e_2 - e_6$ & $x_2 - x_6 \leq -12$\\
 $e_3 - e_6$ & $x_3 - x_6 \leq -10$\\
 $e_4 - e_5$ & $x_4 - x_5 \leq -5$\\
 $e_4 - e_6$ & $x_4 - x_6 \leq -14$\\
 $e_5 - e_6$ & $x_5 - x_6 \leq -5$\\
 $-e_1 + e_5$ & $-x_1 + x_5 \leq 23$\\
 $-e_2 + e_5$ & $-x_2 + x_5 \leq 21$\\
 $-e_2 + e_6$ & $-x_2 + x_6 \leq 28$\\
 $e_2+e_4-e_5$ & $x_2 + x_4-x_5 \leq 8$\\
 $e_2+e_4-e_6$ & $x_2 + x_4 - x_6 \leq 1$\\
 $e_3+e_4-e_6$ & $x_3 + x_4 - x_6 \leq 3$\\
 $e_3+e_5-e_6$ & $x_3 + x_5 - x_6 \leq 13$\\
 $e_1 - e_2 + e_6$ & $x_1 - x_2 + x_6 \leq 33$\\
 $e_1 - e_2 - e_4$ & $x_1 - x_2 - x_4 \leq -7$\\
 $e_1 - e_3 - e_4$ & $x_1 - x_3 - x_4 \leq -12$\\
 $e_1 - e_3 - e_5$ & $x_1 - x_3 - x_5 \leq -19$\\
 $e_2 - e_3 - e_5$ & $x_2 - x_3 - x_5 \leq -17$\\
 $-e_1 + e_5 - e_6$ & $-x_1 + x_5 - x_6 \leq -7$\\
 $e_1+ e_2 - e_3 - e_5$ & $x_1 + x_2 - x_3 - x_5 \leq -9$\\
 $e_2 +e_4 - e_5 - e_6$ & $x_2 + x_4 - x_5 - x_6 \leq -19$\\
 $e_1- e_2 - e_4 + e_6$ & $x_1 - x_2 - x_4 + x_6 \leq 26$\\
 $e_2 - e_3 +e_4 - e_5$ & $x_2 - x_3 + x_4 - x_5 \leq -7$\\
 $-e_1+ e_3 + e_5 -e_6$ & $-x_1 + x_3 + x_5 - x_6 \leq 11$\\
\hline
\end{tabular}
\bigskip
\caption{Rays and Inequalitiesfor $F_{3,7}$}
\label{42}
\end{table}

\section{Cluster Algebras} \label{cluster}

Cluster algebras are commutative algebras endowed with a certain combinatorial structure, introduced in \cite{Cluster1} and expected to be relevant in studying total positivity and homogeneous spaces, such as Grassmannians.

We will not attempt to give a precise definition of a cluster algebra here, but will rather describe their key properties. Slightly varying definitions can be found in \cite{Cluster1}, \cite{Cluster2} and \cite{Scott}; we follow \cite{Scott} but do not believe these small variations are important.

A cluster algebra is an algebra $\AAA$ over a field $k$, which in our examples can be thought of as $\R$. Additionally, a cluster algebra carries two subsets $C$ and $X \subset \AAA$, known as the coefficient variables and the cluster variables. $C$ is finite, but $X$ may be finite or infinite.  If $X$ is finite, $\AAA$ is known as a cluster algebra of finite type. There is also a nonnegative integer $r$ associated to a cluster algebra and known as the rank of the algebra.

There is a pure $(r-1)$-dimensional simplicial complex called the 
{\it cluster complex} whose vertices are the elements of $X$ and 
whose maximal simplices are called {\it clusters}. 
We will denote the cluster complex by $S(\AAA)$. 
If $x \in X$ and $\Delta \in S(\AAA)$ is a cluster containing $x$, 
there is always a unique cluster $\Delta'$ with 
$\Delta \cap \Delta'=\Delta \setminus \{ x \}$. 
Let $\Delta'=(\Delta \setminus \{x\}) \cup \{ x' \}$.  
Then there is a relation $x x'=B$ where $B$ is a binomial in the variables of $(\Delta \cap \Delta') \cup C$.

For any $x \in X$ and any cluster $\Delta$, $x$ is a subtraction-free rational expression in the members of $\Delta \cup C$ and is also a Laurent polynomial in the members of $\Delta \cup C$. 
Conjecturally, this Laurent polynomial has non-negative coefficients.
Note that this conjecture does not follow from the preceding sentence: 
$\frac{x^3+y^3}{x+y}=x^2-xy+y^2$ is a subtraction-free expression in $x$ 
and $y$, and a polynomial in $x$ and $y$, but it 
is not a polynomial with positive coefficients.

It was demonstrated in \cite{Scott} that the coordinate rings of Grassmannians 
have natural cluster algebra structures.
Usually these cluster algebras are of infinite type, making them hard to work with in practice, but in the cases of $Gr_{2,n}$, $Gr_{3,k}$ for $k \leq 8$ and their duals, we get cluster algebras of finite type.

In the case of $Gr_{2,n}$, the coefficient set $C$ is  $\{ \Delta_{12}, \Delta_{23}, \ldots, \Delta_{(n-1) n}, \Delta_{1n} \}$ and the set of cluster variables $X$ is $\{ \Delta_{ij} : i<j \text{ and }i-j \not\congruent \pm 1 \mod n \}$. (Note that these $\Delta$'s are Pl\"ucker coordinates and not simplices.)
Label the vertices of an $n$-gon in clockwise order
with the indices $\{1, 2, \cdots, n \}$ and associate to each member 
of $X \cup C$ the corresponding chord of the $n$-gon.
The  clusters of $Gr_{2,n}$
correspond to the collections of chords which triangulate the $n$-gon. 
Thus, $S(\AAA)$ in this example is (as an abstract simplicial complex) 
isomorphic to the dual of the
associahedron.  Since we have shown
that the fan of $\Trop^+ Gr_{2,n}$ is combinatorially equivalent to 
the normal fan of the associahedron, it follows that
$\Trop^+ Gr_{2,n}$ is (combinatorially) the cone on
$S(\AAA)$. 

In the case of $Gr_{3,6}$, the coefficient set $C$ is equal to
$\{\Delta_{123}, \Delta_{234}, \Delta_{345}, \Delta_{456}, \Delta_{561},\\ \Delta_{612} \}$. $X$ contains the other $14$ Pl\"ucker coordinates, but it also contains two unexpected elements: $\Delta_{134} \Delta_{256} - \Delta_{156} \Delta_{234}$ and $\Delta_{236} \Delta_{145} - \Delta_{234} \Delta_{156}$. 
By definition, all Pl\"ucker coordinates are positive on
the totally positive Grassmannian, so by the results above on 
subtraction-free rational expressions, these new coordinates are positive on the totally positive Grassmannian as well.   

The new coordinates turn out to be Laurent polynomials with positive 
coefficients in the region variables of Section \ref{Parameterize}. Thus, 
we can tropicalize these Laurent polynomials and associate a fan to 
each of them.
When we refine $F_{3,6}$ by these fans,
the refinement 
subdivides the two bipyramids and yields precisely the normal fan to the
$D_4$ associahedron, which is again the cone over $S(\AAA)$.

In the case of $Gr_{3,7}$, $C$ again consists of $\{ \Delta_{i(i+1)(i+2)} \}$ where 
indices are modulo $7$. $X$ contains all of the other Pl\"ucker variables and the 
pullbacks to $Gr_{3,7}$ of the two new cluster variables of $Gr_{3,6}$, 
along with
the $7$ rational coordinate projections $Gr_{3,7} \to Gr_{3,6}$. Thus, $X$ contains $28$ Pl\"ucker variables and $14$ other variables.

As in the case of $Gr_{3,6}$, the 14
new variables are Laurent polynomials with positive coefficients
in the region variables of Section \ref{Parameterize}, so 
to each one we can associate a corresponding fan.  When we refine
$F_{3,7}$ by these 14 new fans, we get a fan combinatorially
equivalent to the fan normal to the $E_6$ associahedron.

We can describe what we have seen in each of these Grassmannian examples in terms of the general language of cluster algebras as follows:

\newtheorem*{ObsBehav}{Observation when $\AAA$ is the Coordinate Ring of a Grassmannian}

\begin{ObsBehav}
Embed $\Spec \AAA$ in affine space by the variables $X \sqcup C$. Then 
$\Trop^+ \Spec \AAA$ is a fan with lineality space of dimension $|C|$. After taking the quotient by this lineality space, we get a simplical fan abstractly isomorphic to the cone over $S(\AAA)$.
\end{ObsBehav}

This observation does not quite hold
for an  arbitrary cluster algebra of finite type. 
For example, if we take the cluster algebra of $Gr_{2,6}$ and set all 
coefficient variables equal to $1$, we get a different cluster algebra
which is still of type $A_3$. However, when we compute
the positive part of the corresponding tropical variety, we get a fan whose 
lineality space has dimension  
$1$, not $0$ as the above would predict. Our fan is 
a cone over a hexagon cross a $1$-dimensional lineality space, 
which is a coarsening of the fan
normal to the type $A_3$ associahedron. 
Based on this and other small examples, it seems that in order to see the 
entire cone over $S(\AAA)$, one needs to use ``enough"
coefficients.

\begin{conjecture}
Let $\AAA$ be a cluster algebra of finite type over $\R$ and $S(\AAA)$
its associated cluster complex.  If the lineality space of $\Trop^+ \Spec \AAA $ 
has 
dimension $|C|$ then $\Trop^+ \Spec \AAA$ 
modulo its lineality space is a simplicial fan abstractly isomorphic to the cone over $S(\AAA)$. If the condition on the lineality space does not hold, 
the resulting fan is a coarsening of the cone over $S(\AAA)$.
\end{conjecture}

\textbf{Remark:} The condition on the lineality space can be restated without mentioning tropicalizations. Consider the torus $(\R^*)^{X \cup C}$ acting on the affine space $\R^{X \cup C}$ and let $G$ be the subgroup taking $\Spec \AAA$ to itself. We want to require that $\dim G=|C|$.

\textbf{Remark:}
In the notation of \cite{Cluster1} and \cite{Cluster2}, the condition on the dimension of the lineality space is equivalent to requiring that the matrix $\tilde{B}$ be of full rank.  We thank Andrei Zelevinsky for pointing this out to
us.

Note how surprising this conjecture is in light of how the two complexes are computed. The fan described in the conjecture is computed as the refinement of a number of fans, indexed by the vertices of $S(\AAA)$. That the rays of this fan, which arise as the intersections of many hypersurfaces, should again be in bijection with the vertices of $S(\AAA)$ is quite unexpected.

We expect an analogous statement to hold for infinite type cluster algebras.

\section{Acknowledgements}
We are grateful to Bernd Sturmfels for suggesting this problem
to us, and for advice throughout our work.   
In addition, we thank Komei Fukuda and Michael Joswig for their help in performing the polyhedral computations and Joshua Scott, Richard Stanley, and 
Andrei Zelevinsky for their useful comments. 

\raggedright

%\end{spacing}
\end{document}